\documentclass[12pt,a4paper]{article}
\usepackage[T1]{fontenc}
\usepackage[latin1]{inputenc}
\usepackage{amsmath,amssymb}
\usepackage[english,frenchb]{babel}

\newcommand{\ddfrac}[2]{\dfrac{\displaystyle #1}{\displaystyle #2}}
\newtheorem{Def}{D{\'e}finition}
\newtheorem{The}{Th{\'e}or{\`e}me}
\newtheorem{Lem}{Lemme}

\begin{document}

 \title{Critical functions on Riemannian manifold wich are locally Euclidien }

\author{ moubinool  omarjee}

\maketitle

\begin{abstract}

In this note we will show the existence of criticals functions for the
functional $I_{f,h}$ where $h$ is a critical function for the metric $g$ 
of $M$ Riemannian manifold of dimension greater than $4$.

\end{abstract}

 INTRODUCTION :
  
Let $M$ a compact Riemannian manifold of dimension greater than $4$ locally Euclidien 
,$h$ a function in  $C^\infty(M)$ , we define for $\phi \in H_1^2(M)-\{0\} $
 $$I_{f,h}(\phi) =
\frac{\int_M {\mid \nabla \phi \mid_g}^2 dv_g +\int_M h{\phi}^2 dv_g
  }{({\int_M f{ \phi}^\frac{2n}{n-2} dv_g})^\frac{2}{2^*}}$$
where $2^*=\frac{2n}{n-2}$ critical Sobolev exponent, $ f \in C^\infty(M,R^*_+)$
we will write

 $$\mu_{f,h} =  \inf\{ \ I_{g,h}(\phi) |  \phi \in
  H_1^2(M) , \phi \not \equiv 0 \}$$

\begin{Def}

We say that $h$ is weakly critical for $g$ if
$$ \mu_{f,h} = \frac{1} {K(n,2)(sup f)^{\frac{n-2}{n}}}$$ 

\end{Def}

\begin{Def}

We say h is critical for $g$ if it is weakly critical
and for any  $k \in C^\infty(M)$ such that
$k\not = h$ we have $ k \leq  h $,
$$\mu_{f,k}< \frac{1} {K(n,2)(sup f)^{\frac{n-2}{n}}}$$

\end{Def}

\begin{Def}

We say that  $\phi \in H_1^2(M)-\{0\} $ is an extremal 
for $h$ if $$I_{f,h}(\phi)=\mu_{  g,h}$$

\end{Def}

Standard elliptic arguments show that an extremal functions belongs to

et qu'elle est soit strictement positive ou strictement
n{\'e}gative, on peut supposer qu'elle est strictement positive.Rechercher
une fonction extr{\'e}male est {\'e}quivalent {\`a} trouver une fonction $\phi \in H_1^2(M)-\{0\} $
telle que 
$$\triangle_g \phi + h\phi =\frac{1} {K(n,2)(sup
f)^{\frac{n-2}{n}}}f\phi ^{2^*-1} \quad et \quad \int_Mf\phi ^{2^*}dv_g=1  $$

\begin{The} 
Soit une vari{\'e}t{\'e} Riemannienne compacte, de dimension plus grande que 4,
localement Euclidienne, c'est {\`a} dire  pour tout $x\in M$, $ Rm_g(x)=0$.
Si $h$ une fonction critique pour $g$ telle que pour tout $x\in M$,
$$h(x)> \frac{(n-2)^2\triangle f(x)}{2(n-1)\omega_n^\frac{2}{n}f(x)}$$ 
Alors il existe une fonction extr{\'e}males pour $h$ .
\end{The}

Preuve:
 Soit $(\alpha_t)_{t\in[0,1]}$ une d{\'e}formation continue de h telle que 

$$\alpha_1 = h$$
$$\forall t \in [0,1[ , \mu_{f,\alpha_t} < 
\frac{1} {K(n,2)(sup f)^{\frac{n-2}{n}}}   \eqno (1) $$  

Prenons $ \alpha_t(x) = th(x) +(1-t) w(x) $ o{\`u}
$w(x)=\frac{(n-2)^2\triangle f(x)}{2(n-1)\omega_n^\frac{2}{n}f(x)}$

les techniques standards elliptiques 

et (1) donne l'existence d'une
famille de fonctions  $ u_t $ de $ C^\infty(M) $, positives 

telles que

$$\forall x \in M ,  \triangle _g u_t(x) + \alpha_t(x) u_t(x) =   \mu_{f,\alpha_t}f(x)u_{t}(x)^{2^* -1}
\eqno(2)$$
$$ \int_M f.u_t^{2^*} dv_g = 1\eqno(3) $$

o{\`u} $\mu_{f,\alpha_t } = \{ \inf\ I_{g,\alpha_t} (\phi) | \phi \in
H_1^2(M)  \phi \not \equiv 0 \}$

La famille $u_t$ est born{\'e}e dans $H_1^2(M)$ 

A extraction pr{\`e}s , il existe une fonction  $ u \geq 0 $ dans  
$H_1^2(M) $ , telle que   

\[{u_t\quad \text {converge vers u }}\quad  \left\{\begin{array}{ll}
                  \text{faiblement} & \mbox{dans} \quad  H_1^2(M)   \\
                  \text{fortement}  & \mbox{dans} \quad  L^2(M)
\end{array}
\right.\]    

   quand $t \rightarrow 1$

$$ \lim_{t\rightarrow 1^-} \mu_{f  ,\alpha_t}
 = \frac{1} {K(n,2)(supf)^{\frac{n-2}{n}}}   \eqno (4) $$

Le principe du maximum donne que la solution $u$ est identiquement nulle 

ou , que $ u > 0$ .Si $u>0$ nous avons trouv{\'e} une fonction extr{\'e}male  

En revanche si $u$ est identiquement nulle, on utilise les techniques
de 

Blow-up .

Introduisons $x_t$ la famille de points en laquelles
$u_t(x_t)=||u_t||_\infty$  

En consid{\'e}rant la famille d'applications $\exp_{x_t}$ qui r{\'e}alise un 

diff{\'e}omorphisme de $B(0,\delta) \subset \bf R^n$ sur la boule
$B(x_t,\delta)$ on va se ramener 

au cas Euclidien, afin d'utiliser l'in{\'e}galit{\'e} de Sobolev Euclidien {\`a}
savoir $$\left(\int_{R^n} |q(x)|^{2^*}dx \right)^\frac{2}{2^*} 
\le K(n,2) \int_{R^n}|\nabla(q(x))|^2 dx$$ 

Et d'aboutir {\`a} la  contradiction $h(x_0)<w(x_0)$ par  des estim{\'e}es
$$ g_t(x) = \exp^*_{x_t}g(x), v_t(x) = u_t(exp_{x_t}(x)) \quad avec\quad
x \in B(0,\delta ) $$

La famille de m{\'e}trique $g_t$ est obtenu comme  le pull-back de la
m{\'e}trique $g$ 

$u_t$ poss{\`e}de un seul point de concentration $x_0$ dont nous rappelons
la 

d{\'e}finition $x_0 \in M$ est un point de concentration pour $u_t$ si
pour tout 

$\delta >0$ 
$$\limsup_{t\rightarrow 1^-}\int_{B(x_0,\delta)}u_t^{2*}dv_g > 0$$ 

Le principe de Moser it{\'e}rative nous dit que les fonctions $u_t$ tend
vers $0$ 

sur le compl{\'e}mentaire d'un voisinage du point de
concentration $x_0$ dans 

la vari{\'e}t{\'e} $M$ .

Soit $\eta$  une fonction $C^{\infty}(B(0,\delta))$ cut-off 

 \[{\eta}= \left\{\begin{array}{ll}
              1 & \mbox{sur} \quad  B(0,\delta/2) \\
              0 & \mbox{sur} \quad   M\backslash B(0,\delta)
              
\end{array}
 \right.\]

Et $|\nabla\eta| \leq C\delta ^{-1}$

On applique la fonction $\eta.v_t$ de support contenu dans
$B(0,\delta)$ dans 

l'in{\'e}galit{\'e} de Sobolev Euclidien 

$$ \left(\int_{R^n} |\eta.v_t|^{2^*}d\xi \right)^\frac{2}{2^*} 
\le K(n,2) \int_{R^n}|\nabla_{\xi}(\eta.v_t)|^2 d\xi
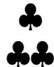
$$ 

\begin{center}
$\clubsuit$\\
$\clubsuit\clubsuit$
\end{center}

\begin{center}
\textbf{ ETAPE 1}
\end{center}

\begin{center}
\textbf{ MONTRONS QUE}
\end{center}

{\large $$\ddfrac{\int_{B(0,\delta)}\eta
  ^2v_t^2\alpha_td\xi}{ \int_{B(0,\delta)}v_t^2d\xi} \leq
\ddfrac{ C^2\delta  ^{-2} \int_{B(0,\delta)\backslash
  B(0,\delta/2)}v_t^2d\xi}
{\int_{B(0,\delta)}v_t^2d\xi} +
 \ddfrac{ A_t(f)}{\int_{B(0,\delta)}v_t^2d\xi} $$

Avec $$A_t(f) =
 \frac{1}{K(n,2)}\left(\frac{1}{||f||_\infty^{\frac{2}{2^*}}}\int_{B(0,\delta)}\eta
 ^2f(\exp_{x_t}(x)) v_t^{2^*}d\xi
- \left(\int_{B(0,\delta)}(\eta
 v_t)^{2^*}d\xi\right)^\frac{2}{2^*} \right)$$\\ }

Un calcul donne 

$$ \int_{B(0,\delta)}|\nabla_{\xi}(\eta.v_t)|^2_{\xi}  d\xi = 
\int_{B(0,\delta)} \partial_{i}(\eta v_t)\partial_{j}(\eta
v_t)d\xi $$
$$=\int_{B(0,\delta)}v_t^2|\nabla \eta|^2_{\xi}d\xi+2\int_{B(0,\delta)}\eta
v_t(\nabla \eta |\nabla v_t)_{\xi} d\xi+\int_{B(0,\delta)}\eta^2|\nabla v_t |^2_{\xi}d\xi $$

En utilisant la fonction cut-off nous obtenons 
$$ \int_{R^n}|\nabla_{\xi}(\eta.v_t)|^2_{\xi}  d\xi \le 
C^2{\delta}^{-2} \int_{B(0,\delta) \backslash B(0,\delta/2)} v_t^2 d\xi +
\int_{B(0,\delta)} \eta^2v_t \triangle_{\xi} v_t d\xi$$

La vari{\'e}t{\'e} est localement Euclidien on a
$\triangle_{\xi}v_t=\triangle_{g_t}v_t$ 
$$
\triangle_{g_t}v_t+\alpha_t(\exp_{x_t}(x))v_t(x)=\mu_{\alpha_t}f(\exp_{x_t}(x))v_t^{2^*-1}(x)$$

On repporte $\triangle_{g_t}v_t$ dans la deuxi{\`e}me int{\'e}grale de
l'in{\'e}galit{\'e} ci-dessus
$$ \int_{R^n}|\nabla_{\xi}(\eta.v_t)|^2_{\xi}  d\xi \le 
C^2{\delta}^{-2} \int_{B(0,\delta) \backslash B(0,\delta/2)} v_t^2 d\xi +
\mu_{\alpha_t} \int_{B(0,\delta)} \eta^2f(\exp_{x_t}(x))v_t^{2^*}d\xi
-$$
$$\int_{B(0,\delta)}\eta^2 \alpha_tv_t^{2}d\xi$$

Cette majoration de $ \int_{R^n}|\nabla_{\xi}(\eta.v_t)|^2_{\xi}
d\xi$ on l'injecte dans 

l'in{\'e}galit{\'e} de Sobolev Euclidien d'o{\`u}
$$ \left(\int_{R^n} |\eta.v_t|^{2^*}d\xi \right)^\frac{2}{2^*} 
\le K(n,2) C^2{\delta}^{-2} \int_{B(0,\delta) \backslash B(0,\delta/2)} v_t^2 d\xi +$$
$$ K(n,2)\mu_{\alpha_t} \int_{B(0,\delta)} \eta^2f(\exp_{x_t}(x))v_t^{2^*}d\xi
-K(n,2) \int_{B(0,\delta)}\eta^2 \alpha_tv_t^{2}d\xi$$

A partir de l{\`a} on a une majoration de $\int_{B(0,\delta)}\eta^2 \alpha_tv_t^{2}d\xi $
$$\int_{B(0,\delta)}\eta^2 \alpha_tv_t^{2}d\xi \le 
C^2{\delta}^{-2} \int_{B(0,\delta) \backslash B(0,\delta/2)} v_t^2 d\xi+
\mu_{\alpha_t} \int_{B(0,\delta)} \eta^2f(\exp_{x_t}(x))v_t^{2^*}d\xi-$$
$$\left(\int_{R^n} |\eta.v_t|^{2^*}d\xi \right)^\frac{2}{2^*} $$

Comme $\mu_{\alpha_t} \le\frac{1} {K(n,2)(sup f)^{\frac{n-2}{n}}} $ on obtient 
$$\ddfrac{\int_{B(0,\delta)}\eta
  ^2v_t^2\alpha_td\xi}{ \int_{B(0,\delta)}v_t^2d\xi} \leq
\ddfrac{ C^2\delta  ^{-2} \int_{B(0,\delta)\backslash
  B(0,\delta/2)}v_t^2d\xi}
{\int_{B(0,\delta)}v_t^2d\xi} +
 \ddfrac{ A_t(f)}{\int_{B(0,\delta)}v_t^2d\xi} $$

Avec  $A_t(f) =
 \frac{1}{K(n,2)}\left(\frac{1}{||f||_\infty^{\frac{2}{2^*}}}\int_{B(0,\delta)}\eta
 ^2f(\exp_{x_t}(x)) v_t^{2^*}d\xi
- \left(\int_{B(0,\delta)}(\eta
 v_t)^{2^*}d\xi\right)^\frac{2}{2^*} \right)$\\

Il nous reste {\`a} {\'e}tudier 

$$\limsup_{t\rightarrow 1^-}\ddfrac{\int_{B(0,\delta)}\eta
  ^2v_t^2\alpha_td\xi}{ \int_{B(0,\delta)}v_t^2d\xi}$$

$$\lim_{t\rightarrow 1} { \ddfrac{ \int_{B(0,\delta)\backslash
  B(0,\delta/2)}v_t^2d\xi}{ \int_{B(0,\delta)}
v_t^2d\xi}} $$

$$\limsup_{t\rightarrow1^-}\ddfrac{ A_t(f)}{\int_{B(0,\delta)}v_t^2d\xi}
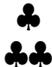
$$

\begin{center}
$\clubsuit$\\
$\clubsuit\clubsuit$
\end{center}

\begin{center}
\textbf{ ETAPE 2}
\end{center}

Nous avons comme dans [2]

$$\lim_{t\rightarrow 1}  
\ddfrac{\int_{B(0,\delta)\backslash B(0,\delta/2)}v_t^2d\xi}{ \int_{B(0,\delta)}v_t^2d\xi}=0 $$

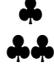
\begin{center}
$\clubsuit$\\
$\clubsuit\clubsuit$
\end{center}

\begin{center}
\textbf{ ETAPE 3}
\end{center}

Etude de 
$$\limsup_{t\rightarrow1^-}\ddfrac{ A_t(f)} {\int_{B(0,\delta)}v_t^2d\xi}$$\\

Posons 
$$D_t(f) = K(n,2)A_t(f)=\frac{1}{||f||^{\frac{2}{2^*}}}\left( \int_{B(0,\delta)}\eta(x)
 ^2f(\exp_{x_t}(x))v_t^{2^*}(x)d\xi \right)-$$
$$ \left(\int_{B(0,\delta)}(\eta v_t)^{2^*}(x)d\xi \right)^\frac{2}{2^*}$$

Posons $\beta_{t}=||u_t||^\frac{2}{2-n} $ tend vers $0$ quand $t$ tend
vers $1^-$.

Pour $x$ appartenant {\`a} la boule $B(0,\frac{\delta}{\beta_t}) \subset
R^n$ ,
$\tilde{g_t}(x)=g_t(\beta_tx)=(exp^*_{x_t}g)(\beta_tx)$ 

$\varphi_t(x)=\beta_t^\frac{n-2}{2}u_t(exp_{x_t}( \beta_t x))$ v{\'e}rifie
l'{\'e}quation 
$$\triangle_{\tilde{g_t}}\varphi_t(x)+\beta_t^2\alpha_t(\exp_{x_t}(\beta_tx))\varphi_t(x)=
\mu_{\alpha_t}f(\exp_{x_t}(\beta_tx))\varphi_t^{2^*-1}(x)$$

Quand $t$ tend vers $1^-$ on a 
$$\triangle_
{\xi}\varphi=\frac{1}{K(n,2)||f||^{\frac{2}{2^*}}}f(\exp_{x_0}(0))\varphi
^{2^* -1}= \frac{f(P)^\frac{2}{n}}{K(n,2)}\varphi^{2^* -1}
 $$ 

o{\`u} $f(\exp_{x_0}(0))=f(x_0)$ .Avec Caffarelli-Gidas-Spruck [4] on
a $$\varphi(x)=\frac{1}{K(n,2)^\frac{n-2}{4}f(P)^\frac{1}{2^*}}\left(1+\frac{|x|^2}{n(n-2)K(n,2)^2}
\right )^{1-\frac{n}{2}}$$

La vari{\'e}t{\'e} est localement Euclidienne la m{\'e}trique $g_t$ et $\xi$
coincident 

localement.

Posons $x=\beta_ty , x^i=\beta_ty^i$ dans $D_t(f)$ on obtient

$$D_t(f)= 
 \frac{1}{||f||^{\frac{2}{2^*}}} \left(
\int_{B(0,\frac{\delta}{\beta_t} )}\eta(\beta_t y )^2
f(\exp_{x_t}(\beta_t y ))v_t^{2^*}(\beta_t y )d\tilde{g_t}  \right)-$$
$$ \left(\int_{B(0,\frac{\delta}{\beta_t} )}( \eta(\beta_t y) v_t(\beta_t y))^{2^*}d\tilde{g_t}
\right)^\frac{2}{2^*} 
 $$

Nous allons appliquer l'in{\'e}galit{\'e} d'Holder {\`a} la premi{\`e}re int{\'e}grale en

{\'e}crivant
$$\eta(\beta_t y
)^2f(\exp_{x_t}(\beta_t y ))v_t^{2^*}(\beta_t y )=
\left( f(\exp_{x_t}(\beta_t y
)) v_t(\beta_t y)^\frac{4}{n-2}\right)
\left(\eta(\beta_t y)^2v_t^2(\beta_t y)\right) $$

D'o{\`u} la majoration de la premi{\`e}re int{\'e}grale intervenant dans $D_t(f)$

$$\frac{1} {||f||^{\frac{2}{2^*}}}\
\int_{B(0,\frac{\delta}{\beta_t} )}\eta(\beta_t y
)^2f(\exp_{x_t}(\beta_t y ))v_t^{2^*}(\beta_t y )d\tilde{g_t}=$$
$$\leq \frac{1}{||f||^{\frac{2}{2^*}}}\left(\int_{B(0,\frac{\delta}{\beta_t} )}\left( f(\exp_{x_t}(\beta_t y
)) v_t(\beta_t y)^\frac{4}{n-2}\right)^\frac{n}{2} d\tilde{g_t}
\right)^\frac{2}{n} \left(\int_{B(0,\frac{\delta}{\beta_t}
  )}\left(\eta(\beta_t y)^2v_t^2(\beta_t y)\right)^\frac{n}{n-2} d\tilde{g_t}\right)^\frac{n-2}{n}$$

En tenant compte du lemme$1$ que nous prouverons {\`a} la fin de l'{\'e}tape 3

\begin{Lem}

Le point de concentration $P$ de la famille $u_t$ coincide avec le 
maximum de $f$, $f(P)=||f||_{\infty}$

\end{Lem}

Nous obtenons la majoration de 
$$D_t(f)\leq  \left(\frac{1}{f(P)^{\frac{2}{2^*}}} \left(\int_{B(0,\frac{\delta}{\beta_t} )} f^\frac{n}{2}(\exp_{x_t}(\beta_t y
)) v_t(\beta_t y)^{2^*} d\tilde{g_t}\right)^\frac{2}{n}  -1\right)\times $$
$$\left(\int_{B(0,\frac{\delta}{\beta_t}
  )} \eta(\beta_t y)^{2^*}v_t^{2^*}(\beta_t y)
d\tilde{g_t}\right)^\frac{2}{2^*}  $$

On {\'e}crit $f^\frac{n}{2}v_t^{2^*}=f(\exp_{x_t}(\beta_t y))
v_t^{2^*}(\beta_t y) f^{\frac{n}{2}-1}(\exp_{x_t}(\beta_t y)) $ afin
de pouvoir 

utiliser ult{\'e}rieurement $\int_M f.u_t^{2^*} dv_g =  1$

On effectue un d{\'e}veloppement de $f^{\frac{n}{2}-1}(\exp_{x_t}(\beta_t y))$ d'o{\`u} 

$$f^\frac{n}{2}v_t^{2^*}=f(\exp_{x_t}(\beta_t y))v_t^{2^*}(\beta_ty)\left(f(x_t)+\partial_if(x_t)\beta_ty^i+
\frac{1}{2}\partial_{ij}f(x_t)\beta ^2_ty^iy^j+o(r^2)\right)^{\frac{n}{2}-1}=$$
$$f(\exp_{x_t}(\beta_t
y))v_t^{2^*}(\beta_ty)(f(x_t))^\frac{n-2}{2}(1+\frac{n-2}{2}\frac{\partial_if(x_t)\beta_ty^i}
{f(x_t)}+\frac{n-2}{4}\frac{\partial_{ij}f(x_t)}{f(x_t)}\beta_t^2y^iy^j+$$
$$\frac{(n-2)(n-4)}{8}(\frac{\partial_if(x_t)\beta_ty^i}{f(x_t)})^2\beta_t^2 y^i{^2}
+o(r^2))$$

Donc l'expression $\left(\int_{B(0,\frac{\delta}{\beta_t} )} f^\frac{n}{2}(\exp_{x_t}(\beta_t y
)) v_t(\beta_t y)^{2^*} d\tilde{g_t}\right)^\frac{2}{n} $ devient 
$$\left(\int_{B(0,\frac{\delta}{\beta_t})}[(f(x_t))^\frac{n-2}{2}f(\exp_{x_t}(\beta_t y))v_t(\beta_t y)^{2^*}+\right.$$
 $$\left.\frac{n-2}{2}(f(x_t))^\frac{n-4}{2}
 \partial_if(x_t)f(\exp_{x_t}(\beta_t y))v_t(\beta_t y)^{2^*}\beta_ty^i+\right.$$
$$\left.\frac{(n-2)(n-4)}{8}(f(x_t))^\frac{n-6}{2}(\partial_if(x_t))^2
f(\exp_{x_t}(\beta_t y))v_t(\beta_t y)^{2^*}\beta_t^2y^i{^2}+\right.$$
$$\left.\frac{n-2}{4}(f(x_t))^\frac{n-4}{2}\partial_{ij}f(x_t)
f(\exp_{x_t}(\beta_t y))v_t(\beta_t
y)^{2^*}\beta_t^2y^iy^j +o(r^2) ] d\tilde{g_t})\right )^\frac{2}{n} $$

Dans l'expression pr{\'e}c{\'e}dente les termes o{\`u} interviennent
$\partial_if(x_t)$ qui tendent 

vers $\nabla f(P)$ vont s'annuler,donc 

$$\limsup_{t\rightarrow
  1^-}\frac{D_t(f)}{\int_{B(0,\delta)}v_t^2d\xi}\leq\limsup_{t\rightarrow
  1^-}\ddfrac{\left(\frac{1}{f(P)^{\frac{2}{2^*}}}\left( \int_{B(0,\frac{\delta}{\beta_t} )} f^\frac{n}{2}
(\exp_{x_t}(\beta_t y)) v_t(\beta_t y)^{2^*} d\tilde{g_t}\right)^\frac{2}{n} -1\right) }
{\int_{B(0,\frac{\delta}{\beta_t})}v_t^2(\beta_ty)d\tilde{g_t} }
\times $$
$$\left(\int_{B(0,\frac{\delta}{\beta_t}
  )} \eta(\beta_t y)^{2^*}v_t^{2^*}(\beta_t y)
d\tilde{g_t}\right)^\frac{2}{2^*} \leq  $$
$$\limsup_{t\rightarrow
  1^-} \ddfrac{\left(\int_{B(0,\frac{\delta}{\beta_t})} \eta(\beta_t y)^{2^*}v_t^{2^*}(\beta_t y)
d\tilde{g_t}\right)^\frac{2}{2^*}}{\int_{B(0,\frac{\delta}{\beta_t})}v_t^2(\beta_t y)d\tilde{g_t} }\times$$
$$(\frac{1}{f(P)^{\frac{2}{2^*}}}\left(((f(x_t))^\frac{n-2}{2} \int_{B(0,\frac{\delta}{\beta_t}
    )}f(\exp_{x_t}(\beta_t y))v_t(\beta_t
    y)^{2^*}d\tilde{g_t} +\right.$$
$$\left.\frac{n-2}{4}(f(x_t))^\frac{n-4}{2}\partial_{ij}f(x_t)\int_{B(0,\frac{\delta}{\beta_t})}
f(\exp_{x_t}(\beta_t y))v_t(\beta_t
y)^{2^*}\beta_t^2y^iy^jd\tilde{g_t}  +\right.$$
$$\left.\partial_if(x_t)F_t^i+o(\delta ^2) \right)^\frac{2}{n} -1)$$

Dans l'expression ci-dessus on majore la premi{\`e}re int{\'e}grale  par
$\int_M f.u_t^{2^*} dv_g $ 

qui vaut $1$
et $\partial_if(x_t) F_t^i$ tend vers 

$0$ quand $t$ tend vers $1^-$
car $\beta_t$ tend vers $0$ quand  $t$ tend vers $1^-$ o{\`u}
$$F_t^i=\frac{n-2}{2}(f(x_t))^\frac{n-4}{2}\int_{B(0,\frac{\delta}{\beta_t})}f(\exp_{x_t}(\beta_t
y))v_t(\beta_ty)^{2^*}\beta_ty^id\tilde{g_t}+$$
$$\frac{(n-2)(n-4)}{8}(f(x_t))^\frac{n-6}{2}\partial_if(x_t)
\int_{B(0,\frac{\delta}{\beta_t})}f(\exp_{x_t}(\beta_t
y))v_t(\beta_ty)^{2^*}\beta_t^2{y^i}^2d\tilde{g_t} $$

nous obtenons 

$$\limsup_{t\rightarrow
  1^-}\frac{D_t(f)}{\int_{B(0,\delta)}v_t^2d\xi} \leq  $$
$$\limsup_{t\rightarrow
  1^-} \ddfrac{\left(\int_{B(0,\frac{\delta}{\beta_t})} \eta(\beta_t y)^{2^*}v_t^{2^*}(\beta_t y)
d\tilde{g_t}\right)^\frac{2}{2^*}}{\int_{B(0,\frac{\delta}{\beta_t})}v_t^2(\beta_t y)d\tilde{g_t} }\times$$
$$(\frac{1}{f(P)^{\frac{2}{2^*}}}\left(((f(x_t))^\frac{n-2}{2} \int_{B(0,\frac{\delta}{\beta_t}
    )}f(\exp_{x_t}(\beta_t y))v_t(\beta_t
    y)^{2^*}d\tilde{g_t} +\right.$$
$$\left.\frac{n-2}{4}(f(x_t))^\frac{n-4}{2}\partial_{ij}f(x_t)\int_{B(0,\frac{\delta}{\beta_t})}
f(\exp_{x_t}(\beta_t y))v_t(\beta_t
y)^{2^*}\beta_t^2y^iy^jd\tilde{g_t}  +\right.$$
$$\left.\partial_if(x_t)F_t^i+o(\delta ^2) \right)^\frac{2}{n} -1)\leq$$
$$\limsup_{t\rightarrow
  1^-} \ddfrac{\left(\int_{B(0,\frac{\delta}{\beta_t})} \eta(\beta_t y)^{2^*}v_t^{2^*}(\beta_t y)
d\tilde{g_t}\right)^\frac{2}{2^*}}{\int_{B(0,\frac{\delta}{\beta_t})}v_t^2(\beta_t y)d\tilde{g_t} }\times$$
$$(\frac{1}{f(P)^{\frac{2}{2^*}}}\left(((f(x_t))^\frac{n-2}{2}  +\right.$$
$$\left.\frac{n-2}{4}(f(x_t))^\frac{n-4}{2}\partial_{ij}f(x_t)\int_{B(0,\frac{\delta}{\beta_t})}
f(\exp_{x_t}(\beta_t y))v_t(\beta_t
y)^{2^*}\beta_t^2y^iy^jd\tilde{g_t}  +\right.$$
$$\left.\partial_if(x_t)F_t^i+o(\delta ^2) \right)^\frac{2}{n} -1)\leq$$
$$\limsup_{t\rightarrow
  1^-} \ddfrac{\left(\int_{B(0,\frac{\delta}{\beta_t})} \eta(\beta_t y)^{2^*}v_t^{2^*}(\beta_t y)
d\tilde{g_t}\right)^\frac{2}{2^*}}{\int_{B(0,\frac{\delta}{\beta_t})}v_t^2(\beta_t y)d\tilde{g_t} }\times$$
$$(\frac{1}{f(P)^{\frac{2}{2^*}}}\left(((f(x_t))^\frac{n-2}{2}  +\right.$$
$$\left.\frac{n-2}{4}(f(x_t))^\frac{n-4}{2}\partial_{ij}f(x_t)\int_{B(0,\frac{\delta}{\beta_t})}
f(\exp_{x_t}(\beta_t y))v_t(\beta_t
y)^{2^*}\beta_t^2y^iy^jd\tilde{g_t}  +\right.$$
$$\left.\partial_if(x_t)F_t^i+o(\delta ^2) \right)^\frac{2}{n} -1)$$

D'o{\`u}
$$\limsup_{t\rightarrow
  1^-}\frac{D_t(f)}{\int_{B(0,\delta)}v_t^2d\xi} \leq $$
$$\limsup_{t\rightarrow
  1^-} \ddfrac{\left(\int_{B(0,\frac{\delta}{\beta_t})} \eta(\beta_t y)^{2^*}v_t^{2^*}(\beta_t y)
d\tilde{g_t}\right)^\frac{2}{2^*}}{\int_{B(0,\frac{\delta}{\beta_t})}v_t^2(\beta_ty)d\tilde{g_t} }\times$$
$$(\frac{f(x_t)^\frac{2}{2^*}}{f(P)^{\frac{2}{2^*}}}
(1+\frac{n-2}{2n}\frac{\partial_{ij}f(x_t)}{f(x_t)}\int_{B(0,\frac{\delta}{\beta_t})}
f(\exp_{x_t}(\beta_ty))v_t(\beta_ty)^{2^*}\beta_t^2y^iy^jd\tilde{g_t}+\frac{2}{n}
\frac{\partial_if(x_t)}{((f(x_t))^\frac{n-2}{2}}
F_t^i+o(\delta ^2))-1) $$

Soit encore 
$$\limsup_{t\rightarrow
  1^-}\frac{D_t(f)}{\int_{B(0,\delta)}v_t^2d\xi} \leq $$
$$\limsup_{t\rightarrow 1^-}
\ddfrac{\frac{f(x_t)^\frac{2}{2^*}}{f(P)^{\frac{2}{2^*}}}\frac{n-2}{2n}\frac{\partial_{ij}f(x_t)}{f(x_t)}
\int_{B(0,\frac{\delta}{\beta_t})}f(\exp_{x_t}(\beta_ty))v_t(\beta_ty)^{2^*}\beta_t^2y^iy^jd\tilde{g_t}}
{\int_{B(0,\frac{\delta}{\beta_t})}v_t(\beta_ty)^2\tilde{g_t}}\times$$
$$(\int_{B(0,\frac{\delta}{\beta_t})} \eta(\beta_t y)^{2^*}v_t^{2^*}(\beta_t y)d\tilde{g_t})^\frac{2}{2^*}$$

On obtient une majoration de la derni{\`e}re int{\'e}grale en introduisant $f$ 
$$(\int_{B(0,\frac{\delta}{\beta_t})} \eta(\beta_t
y)^{2^*}v_t^{2^*}(\beta_t y)d\tilde{g_t})^\frac{2}{2^*}=$$
$$(\int_{B(0,\frac{\delta}{\beta_t})}\frac{1}{f(\exp_{x_t}(\beta_ty)}f(\exp_{x_t}(\beta_ty))  \eta(\beta_t
y)^{2^*}v_t^{2^*}(\beta_t y)d\tilde{g_t})^\frac{2}{2^*}\leq$$
$$(\sup_{B(0,\frac{\delta}{\beta_t})}\frac{1}{f(\exp_{x_t}(\beta_ty)})^\frac{2}{2^*}
(\int_{B(0,\frac{\delta}{\beta_t})}
f(\exp_{x_t}(\beta_ty)\eta(\beta_t
y)^{2^*}v_t^{2^*}(\beta_t y)d\tilde{g_t})^\frac{2}{2^*}\leq$$
$$(\sup_{B(0,\frac{\delta}{\beta_t})}\frac{1}{f(\exp_{x_t}(\beta_ty)})^\frac{2}{2^*}
(\int_{M}f(x)u_t^{2^*}(x)
dg(x))^\frac{2}{2^*}\leq(\sup_{B(0,\frac{\delta}{\beta_t})}\frac{1}{f(\exp_{x_t}(\beta_ty)})^\frac{2}{2^*}$$

Comme 
$v_t(\beta_ty)=u_t(\exp_{x_t}(\beta_ty))=\beta_t^{1-\frac{n}{2}}\varphi_t(y)$,
 $v_t^2(\beta_ty)=\beta_t^{2-n}\varphi_t^2(y)$

$v_t^{2^*}(\beta_ty)=\beta_t^{-n}\varphi_t^{2^*}(y)$ 

Donc 
$$\limsup_{t\rightarrow
  1^-}\frac{D_t(f)}{\int_{B(0,\delta)}v_t^2d\xi} \leq $$
$$
\ddfrac{\frac{n-2}{2n^2}\frac{\triangle f(P)}{f(P)}\int_{R_+}f(P)\varphi^{2^*}(r)r^2r^{n-1}dr}
{\int_{R_+}\varphi^2(r)r^{n-1}dr}\times
(\frac{1}{f(P)})^\frac{2}{2^*}=$$
$$\frac{n-2}{2n^2}\frac{\triangle
  f(P)}{f(P)^\frac{2}{2^*}}\ddfrac{\int_{R_+}f(P)\varphi^{2^*}(r)r^{n+1}dr}
{\int_{R_+}\varphi^2(r)r^{n-1}dr}$$

Car,
$$\limsup_{t\rightarrow 1^- B(0,\frac{\delta}{\beta_t})}\frac{1}{f(\exp_{x_t}(\beta_ty)}=
\limsup_{t\rightarrow 1^-  B(0,\frac{\delta}{\beta})}\frac{1}{f(\exp_{x_t}(0))}=\frac{1}{f(x_0)}=\frac{1}{f(P)}$$

o{\`u} $\varphi(r)=\frac{1}{K(n,2)^\frac{n-2}{4}f(P)^\frac{1}{2^*}}\left(1+\frac{r^2}{n(n-2)K(n,2)}
\right )^{1-\frac{n}{2}}$

Il nous reste {\`a} estimer
$$\ddfrac{\int_{R_+}\varphi^{2^*}(r)r^{n+1}dr}{\int_{R_+}\varphi^2(r)r^{n-1}dr}=$$
$$\frac{1}{K(n,2)f(P)^\frac{2}{n}}\ddfrac{\int_{R_+}\frac{r^{n+1}}{(1+br^2)^n}dr}
{\int_{R_+}\frac{r^{n-1}}{(1+br^2)^{n-2}}dr} \quad avec \quad b=\frac{1}{n(n-2)K(n,2)^2}$$

On int{\`e}gre par parties deux fois le num{\'e}rateur ce qui donne 
$$\int_{R_+}\frac{r^{n+1}}{(1+br^2)^n}dr=\frac{n}{4(n-1)b^2}\int_{R_+}\frac{r^{n-3}}{(1+br^2)^{n-2}}dr$$

D'o{\`u} 
$$\ddfrac{\int_{R_+}\varphi^{2^*}(r)r^{n+1}dr}{\int_{R_+}\varphi^2(r)r^{n-1}dr}\leq\frac{n}{4(n-1)}
\frac{n^2(n-2)^2K(n,2)^4}{K(n,2)f(P)^\frac{2}{n}}$$

Finalement en tenant compte que $K(n,2)^2=\frac{4}{n(n-2)\omega_n^\frac{2}{n}}$
$$\limsup_{t\rightarrow
  1^-}\frac{D_t(f)}{\int_{B(0,\delta)}v_t^2d\xi} \leq $$
$$\frac{n-2}{2n^2}\frac{\triangle
  f(P)}{f(P)^\frac{2}{2^*}}\frac{n}{4(n-1)}\frac{n^2(n-2)^2K(n,2)^4}{K(n,2)f(P)^\frac{2}{n}} =$$
$$\frac{(n-2)^2K(n,2)\triangle f(P)}{2(n-1)\omega_n^\frac{2}{n}f(P)}$$

Avec $A_f(t)$

$$\limsup_{t\rightarrow
  1^-}\frac{A_f(t)}{\int_{B(0,\delta)}v_t^2d\xi} \leq \frac{(n-2)^2\triangle f(P)}{2(n-1)\omega_n^\frac{2}{n}f(P)}$$

Et pour terminer  voici la preuve du lemme$\quad 1$

Par l'absurde si $f(P)<||f||$  , comme $f$ est continue et que $x_t$

tend vers $1^-$ 
on choisit $\delta=\delta_0>0$ petit,et, $t=t_0$ proche de $1^-$
de sorte 

que $$\frac{\sup_{B(x_{t_0},\delta_0)}f(x)}{||f|| }<\frac{1}{3}
\quad {et}\quad
\frac{K(n,2)C^2}{\delta_0}\ddfrac{\int_{B(0,\delta)}v_{t_0}
^2d\xi}{\left(\int_{B(0,\delta)}\eta
  ^{2^*}v_{t_0}
^{2^*}d\xi\right)^\frac{2}{2^*}}<\frac{1}{3}$$

On part de l'{\'e}quation $$\triangle v_t
+\beta_t^2\alpha_t(\exp_{x_t}(x))v_t=\mu_{\alpha_t}f(\exp_{x_t}(x))v_t^{2^*-1}$$
o{\`u} $v_t(x)=u_t(\exp_{x_t}(x))$ ,on multiplie l'{\'e}quation par $\eta
^2v_t$ puis on int{\`e}gre sur la boule $B(0,\delta)$
$$\int_{B(0,\delta)}\eta^2v_t\triangle v_td\xi \leq
\mu_{\alpha_t}\int_{B(0,\delta)}\eta^2f(\exp_{x_t}(x))v_t^{2^*}d\xi
\quad \eqno 1'$$

D'autre part en int{\'e}grant par parties 
$$\int_{B(0,\delta)}\eta^2v_t\triangle v_td\xi=2\int_{B(0,\delta)}\eta
v_t(\nabla\eta|\nabla v_t)d\xi+\int_{B(0,\delta)}\eta^2|\nabla
v_t|^2d\xi$$
Sachant que 
$$\int_{B(0,\delta)}|\nabla(\eta
v_t)|^2d\xi=\int_{B(0,\delta)}v_t^2|\nabla
\eta|^2d\xi+\int_{B(0,\delta)}\eta ^2|\nabla
v_t|^2d\xi+2\int_{B(0,\delta)}\eta v_t(\nabla\eta|\nabla v_t)d\xi$$

Nous arrivons {\`a} 
$$\int_{B(0,\delta)}|\nabla(\eta
v_t)|^2d\xi=\int_{B(0,\delta)}v_t^2|\nabla \eta|^2d\xi
+\int_{B(0,\delta)}\eta^2v_t\triangle v_td\xi \eqno 2'$$

On utilise l'in{\'e}galit{\'e} de Sobolev Euclidien 
$$\left(\int_{B(0,\delta)}\eta
  ^{2^*}v_t^{2^*}d\xi\right)^\frac{2}{2^*}\leq K(n,2)
\int_{B(0,\delta)}|\nabla(\eta v_t)|^2d\xi \eqno 3'$$ 

$1'$  $2'$ $3'$ donnent 

$$\left(\int_{B(0,\delta)}\eta
  ^{2^*}v_t^{2^*}d\xi\right)^\frac{2}{2^*}\leq
K(n,2)\int_{B(0,\delta)}v_t^2|\nabla \eta|^2d\xi+
K(n,2)              \int_{B(0,\delta)}\eta^2v_t\triangle v_td\xi\eqno4'$$
 
Comme $\eta
^2f(\exp_{x_t})v_t^{2^*}=f(\exp_{x_t})v_t^\frac{4}{n-2}.\eta ^2v_t^2$
on applique l'in{\'e}galit{\'e} d'Holder avec
$p=\frac{2^*}{2}$,$q=\frac{n}{2}$ d'o{\`u}

$$\int_{B(0,\delta)}\eta^2f(\exp_{x_t})v_t^{2^*}d\xi\leq
\left(\int_{B(0,\delta)}\eta
  ^{2^*}v_t^{2^*}d\xi\right)^\frac{2}{2^*}\left(\int_{B(0,\delta)}
f^\frac{n}{2}(\exp_{x_t})v_t^{2^*}d\xi\right)^\frac{2}{n} \eqno5'$$

$1'$ $4'$ $5'$ donnent 

$$\left(\int_{B(0,\delta)}\eta
  ^{2^*}v_t^{2^*}d\xi\right)^\frac{2}{2^*}\leq
K(n,2)\int_{B(0,\delta)}v_t^2|\nabla \eta|^2d\xi+$$
$$
K(n,2)\mu_{\alpha_t}\left(\int_{B(0,\delta)}\eta
  ^{2^*}v_t^{2^*}d\xi\right)^\frac{2}{2^*}\left(\int_{B(0,\delta)}
f^\frac{n}{2}(\exp_{x_t})v_t^{2^*}d\xi\right)^\frac{2}{n}\eqno6'$$

Comme $|\nabla \eta| \leq \frac{C}{\delta}$ et $\mu_{\alpha_t}\le
\frac{1}{K(n,2)||f||^\frac{n-2}{n}}$ et

$$\left(\int_{B(0,\delta)}\eta
  ^{2^*}v_t^{2^*}d\xi\right)^\frac{2}{2^*}\leq \frac{K(n,2)C^2}{\delta
  ^2}\int_{B(0,\delta)}v_t^2d\xi                                            +$$
$$\frac{1}{||f||^\frac{n-2}{n}}\left(\int_{B(0,\delta)}
f^\frac{n}{2}(\exp_{x_t})v_t^{2^*}d\xi\right)^\frac{2}{n}\left(\int_{B(0,\delta)}\eta
  ^{2^*}v_t^{2^*}d\xi\right)^\frac{2}{2^*} \eqno7'$$

En utilisant 
 $f^\frac{n}{2}v_t^{2^*}=f^\frac{n-2}{2}fv_t^{2^*}$ on a 
$$\left(\int_{B(0,\delta)}
f^\frac{n}{2}(\exp_{x_t})v_t^{2^*}d\xi\right)^\frac{2}{n}\leq
\left(\sup_{B(x_t,\delta)}f(x)^\frac{n-2}{n}\right)^
\frac{2}{n}\left(\int_{B(0,\delta)}f(\exp_{x_t}(x))
v_t^{2^*}d\xi\right)^\frac{2}{n}$$
La derni{\`e}re int{\'e}grale est major{\'e}e par $\int_Mfu_t^{2^*}=1$

$7'$ devient 

$$\left(\int_{B(0,\delta)}\eta
  ^{2^*}v_t^{2^*}d\xi\right)^\frac{2}{2^*}\leq \frac{K(n,2)C^2}{\delta
  ^2}\int_{B(0,\delta)}v_t^2d\xi                                            +$$
$$\frac{ \sup_{B(x_t,\delta)}f(x)^\frac{n-2}{n} }{||f||^\frac{n-2}{n}}                                                                                 \left(\int_{B(0,\delta)}\eta
  ^{2^*}v_t^{2^*}d\xi\right)^\frac{2}{2^*} \eqno8'$$

Soit encore 

$$1\leq \frac{K(n,2)C^2}{\delta ^2}\ddfrac{\int_{B(0,\delta)}v_t^2d\xi}{\left(\int_{B(0,\delta)}\eta
  ^{2^*}v_t^{2^*}d\xi\right)^\frac{2}{2^*}}+\left(\frac{\sup_{B(x_t,\delta)}f(x)}{||f||}\right)^\frac{n-2}{n}$$
En choisissant $\delta=\delta_0$ et $t=t_0$ dans la pr{\'e}c{\'e}dente in{\'e}galit{\'e} nous avons la
contradiction $1\leq \frac{1}{3} + \frac{1}{3}$ 

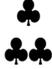
\begin{center}
$\clubsuit$\\
$\clubsuit\clubsuit$
\end{center}

\begin{center}
\textbf{ ETAPE 4}
\end{center}

Arriver {\`a}  une contradiction

En faisant tendre $t$ vers $1$ dans l'in{\'e}galit{\'e} suivante
  $$\ddfrac{\int_{B(0,\delta)}\eta
  ^2v_t^2\alpha_td\xi}{\int_{B(0,\delta)}v_t^2d\xi} \leq
\ddfrac{C^2\delta  ^{-2} \int_{B(0,\delta)\backslash
  B(0,\delta/2)}v_t^2d\xi}
{\int_{B(0,\delta)}v_t^2d\xi} +
 \ddfrac{ A_t(f)}{\int_{B(0,\delta)}v_t^2d\xi} $$

$$h(P)\leq \frac{(n-2)^2\triangle f(P)}{2(n-1)\omega_n^\frac{2}{n}f(P)} $$


I thanks M.Vaugon for suggesting  me this problem.

\end{document}